\documentclass[12pt]{article}
\usepackage{amssymb}

\usepackage[english]{babel}
\usepackage{latexsym}
\usepackage{amsfonts}
\usepackage{amsmath}
\usepackage{amscd}

\newcommand{\FAT}[1]{\mbox{{$\mathbb{#1}$}}}

\newcommand{\CC}{\FAT{C}}

\newcommand{\ZZ}{\FAT{Z}}

\newcommand{\Hom}{{\rm Hom}}

\newcommand{\im}{{\rm Im}}
\newcommand{\Ker}{{\rm Ker}}

\newcommand{\Aut}{\mathop{\rm Aut}\nolimits}
\newcommand{\Inn}{\mathop{\rm Inn}\nolimits}
\newcommand{\Irr}{\mathop{\rm Irr}\nolimits}
\newcommand{\qed}{\nolinebreak\hfill\rule{2mm}{2mm}%
\par\medbreak}

\newtheorem{thm}{Theorem}

\newtheorem{lemma}[thm]{Lemma}

\newtheorem{cor}[thm]{Corollary}

\title{Enumeration of Lifts of Commuting Elements of a Group}
\author{Michael Natapov \\
\small Department of Mathematics, Indiana University, Bloomington,
IN 47405 \\
\small mnatapov@indiana.edu\\
\\
 \small and \\
\\
Vladimir Turaev \\
\small Department of Mathematics, Indiana University, Bloomington,
IN 47405 \\
\small vturaev@yahoo.com}

\date{}                                         

\begin{document}

\maketitle

\section{Introduction}

Frobenius \cite{Fr} computed the number of
pairs of commuting elements of a finite group $\Gamma$ to be  $
c\,\vert \Gamma\vert  $, where $c$ is the number of conjugacy
classes in $\Gamma$ and   the vertical bars denote the cardinality
of a set. This formula  was generalized by  Mednykh \cite{Me} who
computed the number of homomorphisms from the fundamental group
$\pi=\pi_1(W)$ of an arbitrary
closed connected oriented surface $W$ to $\Gamma$:

\begin{equation}\label{FroMed}\vert  {\Hom} (\pi , \Gamma)\vert\, =\, \vert  \Gamma \vert  \,
\sum_{\rho\in {\rm Irr} (\Gamma)} \,\left (\frac{\vert \Gamma
\vert}{\dim \, \rho}\right )^{2d-2}\, .
\end{equation}
Here ${\rm Irr} (\Gamma)$ is the set of equivalence classes of
irreducible complex linear representations of $\Gamma$, and $d$ is
the genus of $W$. For $W=S^1\times S^1$, this gives  the Frobenius
formula   since   $ \pi_1(S^1\times S^1)=\ZZ^2$ and  $\vert {\rm
Irr} (\Gamma)\vert=c$. Formula (\ref{FroMed}) can be proved by
algebraic means, see Jones \cite{Jo}, or deduced from Topological
Quantum Field Theory in dimension 2, see Freed and Quinn \cite{FQ}.

The Frobenius-Mednykh formula (\ref{FroMed}) was generalized by
Turaev \cite{T} as follows. Consider a  group epimorphism $q:G'\to
G$ with finite kernel $\Gamma$. Fix a homomorphism
$g:\pi=\pi_1(W)\to G$. Let ${\rm Hom}_g(\pi, G')$ be the set of all
lifts of $g$ to $G'$, that is the set of homomorphisms $g':\pi\to
G'$ such that $q g'=g$. Since   $\pi$ is finitely generated, the set
${\rm Hom}_g(\pi, G')$ is finite.     Note that the action of $G$
  on $\Gamma$ by outer automorphisms induces an action
of $G$ on ${\rm Irr} (\Gamma)$. The stabilizer of $\rho \in {\rm
Irr} (\Gamma)$ under this action is denoted $G_\rho$. Then
\begin{equation}\label{Turaev1} |{\Hom}_g(\pi, G')|=\,
\vert  \Gamma \vert  \, \sum_{\rho\in {\rm Irr} (\Gamma), G_{\rho} \supset \im(g)} \,\left
(\frac{\vert \Gamma \vert}{\dim \, \rho}\right )^{2d-2} \,g^*(\zeta_\rho) ([W])\, ,
\end{equation}
where $\zeta_\rho\in H^2(G_\rho; \CC^\times) $ is a   cohomology
class introduced in \cite{T}, and $g^*(\zeta_\rho) ([W])\in
\CC^\times$ is the evaluation of $g^*(\zeta_\rho) \in H^2(\pi;
\CC^\times)$ on the fundamental class $[W]\in H_2(W{})=H_2(\pi{})$
of $W$.  Here and below the group of coefficients in homology is
$\ZZ$. The evaluation in question  is induced by the bilinear form
$\CC^\times \times \ZZ\to \CC^\times, (z,n)\mapsto z^n$. For $G=1,
G'=\Gamma$, we recover Formula (\ref{FroMed}).

For $W=S^1\times S^1$, Formula (\ref{Turaev1}) computes the number
of pairs of commuting elements of $G'$ projecting  to a given pair
of commuting elements of $G$. More precisely, let $\alpha, \beta \in
G$ be such that $\alpha \beta = \beta \alpha$. Set
$$L(\alpha, \beta, q) = \{a,b \in G'\, | \, q(a) = \alpha, q(b) = \beta,  \mbox{ and } ab=ba\}\, .$$
Clearly,
 $L(\alpha, \beta, q) = {\Hom}_g(\ZZ^2, G') $,
where $g:\ZZ^2=\pi_1(S^1\times S^1)\longrightarrow G$ carries the
  generators of $ \ZZ^2$ to $\alpha$ and $\beta$,
respectively. Formula (\ref{Turaev1})   implies that
\begin{equation}\label{Turaev1-} |L(\alpha, \beta, q)|=\,
\vert  \Gamma \vert  \, \sum_{\rho\in {\rm Irr} (\Gamma), G_{\rho} \supset \{\alpha, \beta\}} \,
 \,g^*(\zeta_\rho) ([W])\, .
\end{equation}
We  use  this equality to compute   $|L(\alpha, \beta, q)|$
 for several small groups $\Gamma$. We focus on the case where
  $\Gamma=Q_8 = \{\pm 1, \pm i, \pm j, \pm k\}$ is the quaternion
group.

\begin{thm}\label{Theorem_Q_8(1)}
Let $1 \rightarrow Q_8 \rightarrow G' \stackrel{q} \rightarrow G
\rightarrow 1$ be a short exact sequence of groups. For any
commuting elements  $\alpha, \beta $ of $ G$,  the set $L(\alpha,
\beta, q)  $   has $0$, $8$, $16$, $24$, or $40$ elements.
\end{thm}

\begin{cor} Let $G'$ be a group with normal subgroup $Q_8 \triangleleft G'$.
Let $a,b $ be commuting elements of $ G'$. Then the number of pairs
$\gamma_1, \gamma_2 \in Q_8$ such that $a\gamma_1$ commutes with $
b\gamma_2   $ is equal to $8$, $16$, $24$, or $40$.
\end{cor}

\begin{thm}\label{Theorem_Q_8(2)}
For each $n \in \{0,8,16,24,   40\}$, there is a short exact
sequence of the form $Q_8 \hookrightarrow G'
\stackrel{q}{\twoheadrightarrow} G$ and commuting elements $\alpha,
\beta $ of $ G$ such that $|L(\alpha, \beta, q)| = n$.
\end{thm}

The existence of commuting lifts of $\alpha, \beta\in G$ to $G'$ may
be approached from a   homological viewpoint. Consider for
simplicity the case where $\alpha, \beta$ generate $G$. Let
  $g:\ZZ^2\to G$ be the homomorphism carrying the
  generators of $ \ZZ^2$ to $\alpha$ and $\beta$,
respectively. Let  $\kappa_{\alpha, \beta} \in H_2(G{})$ be the
image of a generator of $H_2(\ZZ^2{})=\ZZ$ under the induced
homomorphism $g_*: H_2(\ZZ^2{}) \to H_2(G{})$. The homology class
$\kappa_{\alpha, \beta}$ is well defined, at least up to
multiplication by $-1$. It is clear that if $\alpha, \beta $ lift to
commuting elements of   $G'$, then $\kappa_{\alpha, \beta}$ lies in
the image of the homomorphism $q_*: H_2(G'{})\to H_2(G{})$. This
yields a homological obstruction to the existence of such a lift. It
is easy to give examples showing that this obstruction may be
non-trivial. The following theorem shows  that, generally speaking,
this obstruction is insufficient.

\begin{thm}\label{Theorem_Q_8(2+)}
There is a short exact sequence of the form $Q_8 \hookrightarrow G'
\stackrel{q}{\twoheadrightarrow} G$ and commuting generators
$\alpha, \beta $ of $ G$ such that $\kappa_{\alpha, \beta}$ lies in
the image of the homomorphism $q_*: H_2(G'{})\to H_2(G{})$ and $
L(\alpha, \beta, q)  = \emptyset$.
\end{thm}

Applying Formula (\ref{Turaev1}) to  surfaces of   genus $\geq 2$,
we obtain the following theorem.

\begin{thm}\label{ed} Let $1 \rightarrow Q_8 \rightarrow G' \stackrel{q} \rightarrow
G \rightarrow 1$ be a short exact sequence of groups. Let $d\geq 2$
and $\alpha_i, \beta_i \in G$, $1 \leq i \leq d$ be such that
$\prod_{i=1}^d [\alpha_i, \beta_i] = 1$. Then there is a family
$\{a_i, b_i \}_{i=1}^d$ of elements of $G'$ such that $q(a_i) =
\alpha_i, q(b_i) = \beta_i$ for all $i$ and $ \prod_{i=1}^d [a_i,
b_i] = 1$. The number of such families is equal to $8^{2d-1}(N \pm
2^{2-2d})$, where $N\in \{1,2,4\}$.
\end{thm}

It would be interesting to find out what   numbers  $8^{2d-1}(N \pm
2^{2-2d})$ with  $N\in \{1,2,4\}$ are realizable as the number of
families as in this theorem.

The work of the second named author was partially supported by the
NSF grant DMS-0707078.

\section{Representations and cohomology classes}\label{section2}

We define in this section the cohomology classes $\zeta_\rho$ used
in the formulas above.

Let $q:G'\to G$ be a group epimorphism with kernel $\Gamma$ (not
necessarily finite). The group $G$ acts on ${\rm Irr} (\Gamma)$ as
follows. For each $\alpha \in G$, choose $ \widetilde{\alpha} \in
q^{-1}(\alpha)\subset G'$. Let $\rho : \Gamma \to GL_n(\CC)$ be an
irreducible representation of degree $n\geq 1$. For $\alpha \in G$,
the map $    \Gamma \to GL_n(\CC), \gamma \mapsto
\rho(\widetilde{\alpha} \gamma \widetilde{\alpha}^{-1}) $ is an
irreducible representation of $\Gamma$ denoted $\rho \alpha$. The
equivalence class of $\rho\alpha$ does not depend on the choice of $
\widetilde{\alpha} $. This defines a right action of $G$ on ${\rm
Irr} (\Gamma)$.

Given an irreducible  representation $\rho:\Gamma\to GL_n(\CC)$, let
$ G _\rho = \{ \alpha \in G   | \, \rho \alpha \sim \rho \} $ be the
stabilizer of its equivalence class. We now define $\zeta_\rho \in
H^2(G_\rho; \CC^\times)$ following \cite{T}. For each $\alpha \in
G_\rho $, there is a matrix $M_\alpha\in GL_n(\CC)$ such that
\begin{equation} \label{bh}
\rho(\widetilde{\alpha} \gamma \widetilde{\alpha}^{-1}) = M_\alpha  \rho (\gamma)  M_\alpha^{-1}
\end{equation}
for all $\gamma \in \Gamma$. The irreducibility of $\rho$ implies
that $M_\alpha$ is unique up to multiplication by a non-zero scalar.
Fix $M_\alpha$ for all $\alpha\in G_\rho$. For any $\alpha, \beta\in
G_\rho$, there a unique $\zeta_\rho(\alpha, \beta) \in \CC^\times$
such that
\begin{equation} \label{cohomology class_1}
  M_{\alpha
\beta} \, \rho (
 (\widetilde{\alpha\beta})^{-1} \widetilde \alpha
\widetilde \beta) = \zeta_\rho(\alpha,
\beta) \, M_\alpha\, M_\beta.
\end{equation}
One checks that $\zeta_\rho:G_\rho \times G_\rho \to \CC^\times$ is
a 2-cocycle. Its cohomology class in  $  H^2(G_\rho; \CC^\times)$
depends only on the equivalence class of $\rho$ and does not depend
on the choice of the conjugating matrices $\{M_\alpha\}_\alpha$. By
abuse of notation, we  denote this cohomology class by the same
symbol $\zeta_\rho$. The arguments below in this section show that
the cohomology class  $q^*(\zeta_\rho)\in H^2(q^{-1}(G_\rho);
\CC^\times)$ has finite order. If   $\Gamma$ is finite, then we can
deduce that $\zeta_\rho$ has finite order in $  H^2(G_\rho;
\CC^\times)$ so that its evaluation on any 2-dimensional homology
class of $G_\rho$ is a root of unity.

If $n=1$, then we can take $M_\alpha=1\in \CC$ for all $\alpha\in G$
and obtain  $\zeta_\rho(\alpha, \beta)=\rho (
 (\widetilde{\alpha\beta})^{-1} \widetilde \alpha
\widetilde \beta)$ for all $\alpha, \beta\in G$. In particular, if
$\rho$ is the trivial 1-dimensional representation of $\Gamma$, then
$\zeta_\rho=1$.

A similar construction   produces cohomology classes of certain
subgroups of the group ${\mathcal A}={\Aut}(\Gamma)$  of
automorphisms of $\Gamma$. We define a right action of ${\mathcal
A}$ on $\Irr(\Gamma)$ by $\rho \varphi   = \rho \circ \varphi $ for
$\varphi \in {\mathcal A}$ and $\rho\in \Irr(\Gamma)$. For an
irreducible representation $\rho:\Gamma\to GL_n(\CC)$, the
stabilizer ${\mathcal A}_\rho\subset {\mathcal A} $ of the
equivalence class of $\rho$ consists of all $\varphi \in {\mathcal
A} $ such that there is ${\mathcal M}_\varphi \in GL_n(\CC)$
satisfying
\begin{equation} \label{matrix}
\rho\varphi(\gamma) = {\mathcal M}_\varphi \, \rho (\gamma)\,
{\mathcal M}_\varphi^{-1}\end{equation} for all $\gamma \in \Gamma$.
The matrix  ${\mathcal M}_\varphi$ is unique up to multiplication by
an element of $\CC^\times$. We fix ${\mathcal M}_\varphi$ for all
$\varphi \in {\mathcal A}_\rho$ and define a 2-cocycle
$\eta_\rho:{\mathcal A}_\rho \times {\mathcal A}_\rho \to
\CC^\times$ by
\begin{equation} \label{cohomology class}
{\mathcal M}_{\varphi  \psi} = \eta_\rho(\varphi, \psi) \, {\mathcal M}_\varphi\, {\mathcal M}_\psi\, ,
\end{equation}
where $\varphi, \psi \in {\mathcal A}_\rho$. The   class of this
cocycle in $ H^2({\mathcal A}_\rho; \CC^\times)$     depends only on
the equivalence class of $\rho$ and does not depend on the choice of
the conjugating matrices $\{{\mathcal M}_\varphi\}_\varphi$. This
cohomology class is also denoted $\eta_\rho$. Taking the determinant
on both sides of (\ref{cohomology class}), we obtain that
$\eta_\rho^n=1$, where $n=\dim\, \rho$. The constructions of
$\zeta_\rho$ and $\eta_\rho$ are related   by the following lemma.

\begin{lemma}\label{le1} Let $q:G'\to G$ be a group epimorphism with kernel $\Gamma$.
  Let $\omega: G' \to {\mathcal A}={\Aut}(\Gamma) $ be  the
homomorphism carrying $a\in G'$ to the automorphism $\gamma\mapsto
a\gamma a^{-1}$ of $\Gamma$. For any irreducible representation
$\rho:\Gamma \to GL_n(\CC)$ of $\Gamma$, we have
$\omega^{-1}({\mathcal A}_\rho)=q^{-1} (G_\rho)$. If
$\omega(G')\subset  {\mathcal A}_\rho$ (or equivalently, if
$G_\rho=G$), then
\begin{equation}\label{Turaev_2++}
q^*(\zeta_\rho)={\omega}^*(\eta_\rho)\, ,
\end{equation}
where $q^* : H^2(G;\CC^\times) \to H^2(G';\CC^\times)$ and
${\omega}^* : H^2({\mathcal A}_\rho;\CC^\times) \to
H^2(G';\CC^\times)$ are the homomorphisms induced by $q$ and
$\omega$, respectively.
\end{lemma}

{\it Proof.} The equality $\omega^{-1}({\mathcal A}_\rho)=q^{-1}
(G_\rho)$ follows from the definitions. For each $\alpha \in G$ fix
$ \widetilde{\alpha} \in q^{-1}(\alpha)\subset G'$  and set
$M_\alpha= \mathcal M_{\omega(\widetilde {\alpha})}$. Formula
(\ref{matrix}) implies    (\ref{bh}).

For   $a\in G'$, set $\gamma_a=
 (\widetilde {q(a)})^{-1} a \in \Gamma$ and $M^+_a=M_{q(a)}
\rho(\gamma_a)\in GL_n(\CC)$. Given $\gamma\in \Gamma$, we deduce
from (\ref{bh}) that
$$\rho(a\gamma a^{-1})=\rho (\widetilde {q(a)} \gamma_a \gamma \gamma_a^{-1} (\widetilde {q(a)})^{-1})=
 M^+_a\, \rho(\gamma)\, (M^+_a)^{-1} \, .$$ Comparing
with   (\ref{matrix})  we obtain  that $M^+_a=k_a \, \mathcal
M_{\omega(a)} $ for some $k_a\in \CC^\times$.

The 2-cocycle $\{\eta_\rho(\omega(a), \omega(b))\}_{a,b\in G'}$ can
be computed from the matrices $\{\mathcal M_{\omega(a)}\}_a$ via the
identity $\mathcal M_{\omega(ab)}=
 \eta_\rho(\omega(a), \omega(b)) \, \mathcal M_{\omega(a)} \, \mathcal M_{\omega(b)}$ for   $a,b\in G'$. The 2-cocycle $\{\zeta_\rho({q(a)},
{q(b)})\}_{a,b\in G'}$ derived from the matrices
$\{M_\alpha\}_{\alpha\in G}$  can be   computed from  $\{M^+_a\}_a$.
Indeed, for $a,b\in G'$, we have
$$\gamma_{ab}= ({\widetilde {q(ab)}})^{-1} \, {\widetilde {q(a)}}\, \gamma_a \, {\widetilde {q(b)}}\, \gamma_b\,  $$
and
$$M^+_{ab}=M_{q(ab)}\, \rho(\gamma_{ab})=
M_{q(a)q(b)}\,
\rho(({\widetilde {q(ab)}})^{-1} \, {\widetilde {q(a)}}\, \gamma_a \, {\widetilde {q(b)}}\, \gamma_b)$$
$$= M_{q(a)q(b)}\,
\rho(({\widetilde {q(ab)}})^{-1} \, {\widetilde {q(a)}}\,   {\widetilde {q(b)}})\,
 \rho( {\widetilde {q(b)}}^{-1} \gamma_a  {\widetilde {q(b)}})\, \rho(
 \gamma_b)$$
$$=
 \zeta_\rho({q(a)},
{q(b)}) \, M_{q(a)}\, M_{q(b)} M_{q(b)}^{-1} \,\rho(\gamma_a) M_{q(b)}\, \rho (\gamma_b)=
 \zeta_\rho({q(a)},
{q(b)}) M^+_a M^+_b\, .$$
 We conclude that the 2-cocycles $\{\zeta_\rho({q(a)},
{q(b)})\}_{a,b\in G'}$ and $\{\eta_\rho(\omega(a),
\omega(b))\}_{a,b\in G'}$ differ by a coboundary.\qed

Formula (\ref{Turaev_2++}) allows us to rewrite Formulas
(\ref{Turaev1}) and (\ref{Turaev1-}) as follows. Suppose that
$g(\pi)=G$ and $g_*([W]) = q_*(\Delta)$ for some $\Delta \in
H_2(G'{})$, where $q_*: H_2(G'{}) \to H_2(G{})$ and $g_*: H_2(\pi{})
\to H_2(G{})$ are the homomorphisms induced by $q:G'\to G$ and
$g:\pi=\pi_1(W)\to G$, respectively. Then we can replace
$g^*(\zeta_\rho) ([W])$ in (\ref{Turaev1}) and (\ref{Turaev1-}) by
${\omega}^*(\eta_\rho) (\Delta)$. Indeed, $$g^*(\zeta_\rho) ([W])=
\zeta_\rho (g_*( [W]))=\zeta_\rho (q_*(\Delta))=q^*(\zeta_\rho)
(\Delta)={\omega}^*(\eta_\rho) (\Delta)\, .$$ In particular, if
$\Gamma$ is finite and $G$ is generated by commuting elements
$\alpha, \beta$, then (\ref{Turaev1-}) gives
\begin{equation}\label{Turaev 2} |L(\alpha, \beta, q)|=\,
\vert  \Gamma \vert  \, \sum_{\rho\in {\rm Irr} (\Gamma), \,{\mathcal A}_\rho \supset  \omega(G')} \,
 \,{\omega}^*(\eta_\rho) (\Delta) \, .
\end{equation}
Note that  $\Delta$ as above necessarily exists if $L(\alpha, \beta,
q)\neq \emptyset$ because in this case  $g$ lifts to a homomorphism
$g': \pi\to G'$ and  we may take $\Delta =g'_*([W])\in H_2(G'{})$.
Note also that if $\dim \rho=1$, then $\eta_\rho=1$ and
${\omega}^*(\eta_\rho) (\Delta)=1$ for any $\Delta  \in H_2(G'{})$.
If $\dim \rho=2$, then   ${\omega}^*(\eta_\rho) (\Delta)=\pm 1$ for
any $\Delta  \in H_2(G'{})$.

\section{Proof of Theorems 1 and 5}

{\it Proof of Theorem \ref{Theorem_Q_8(1)}.} It is well known that
in the generators $a=i$, $b=j$,  the group  $Q_8 $ has the
presentation $\{a,b\, | \, a^4=1, a^2=b^2, bab^{-1}=a^{-1}\}$.
Therefore the abelianization of $Q_8$ gives $\ZZ_2 \times \ZZ_2$.
Consequently,   $Q_8$ has four homomorphisms to
$\CC^\times=GL_1(\CC)$, one trivial $\chi_0$ and three non-trivial
$\chi_1, \chi_2, \chi_3$. The equality $\sum_{\rho\in \Irr Q_8}
(\dim \rho)^2 =8$ implies that besides the 1-dimensional
representations, the set $\Irr Q_8$ has only one element which is
the equivalence class of an irreducible 2-dimensional representation
$\rho_0$. This representation can be described explicitly  by
$$i\mapsto \left(
             \begin{array}{cc}
               i & 0 \\
               0 & -i \\
             \end{array}
           \right), \quad j\mapsto \left(
             \begin{array}{cc}
               0 & 1 \\
               -1 & 0 \\
             \end{array}
           \right), \quad \mbox{and\ } \,\,\, k=ij\mapsto \left(
             \begin{array}{cc}
               0 & i \\
               i & 0 \\
             \end{array}
           \right).
$$

Consider the group  ${\mathcal A}=\Aut Q_8$ and its subgroups
$\{{\mathcal
 A}_{\chi_i}\}_{i=0}^3$ and ${\mathcal A}_{\rho_0}$. Clearly,  $  {\mathcal
A}_{\chi_0}= {\mathcal A}_{\rho_0}={\mathcal A}$. The obvious
equality $\chi_3=\chi_1  \chi_2$ implies that the intersection of
any two of the sets $\{{\mathcal A}_{\chi_i}\}_{i=1,2,3}$ is equal
to the intersection of all three of these sets.

Replacing if necessary $G$ by its subgroup generated by $\alpha,
\beta$ and replacing $G'$ by the pre-image of this subgroup, we can
assume that $\alpha$ and $\beta$ generate $G$. Let    $\omega: G'
\to {\mathcal A}=\Aut Q_8$ be the homomorphism defined in Lemma
\ref{le1}.   The remarks at the end of Section \ref{section2} show
that ${\omega}^*(\eta_{\chi_i}) (\Delta)=1$ for $i=0,1,2,3$ and
${\omega}^*(\eta_{\rho_0}) (\Delta)=\pm 1$ for any $\Delta\in
H_2(G'{})$.   If   $L(\alpha, \beta, q)= \emptyset$, then
$|L(\alpha, \beta, q)|=0$ and we are done.   If $L(\alpha, \beta,
q)\neq \emptyset$, then by (\ref{Turaev 2}),

\begin{equation} \label{formula_Q_8}
|L(\alpha, \beta, q)| = 8( 1 + \left\{
                                                                             \begin{array}{ll}
                                                                               0, & {\omega(G') \nsubseteq {\mathcal A}_{\chi_i} \mbox{ \rm for }i=1,2,3;} \\
                                                                               1, & {\omega(G') \subseteq {\mathcal A}_{\chi_i} \mbox{ \rm  for only one }i\in \{1,2,3\};} \\
                                                                               3, & {\omega(G') \subseteq {\mathcal A}_{\chi_i}  \mbox{ \rm for }i=1,2,3.}
                                                                             \end{array}
                                                                           \right\} \pm 1
).
\end{equation}
The possible values for the right-hand side are $0,8,16,24$, and
$40$.\qed

{\it Proof of Theorem \ref{ed}.} Replacing if necessary $G$ by its
subgroup generated by $\{\alpha_i, \beta_i\}_{i=1}^d$ and replacing
$G'$ by the pre-image of this subgroup, we can assume that the set
$\{\alpha_i, \beta_i\}_{i=1}^d$ generates $G$.  Let  $W$ be a closed
connected oriented surface of genus $d$. Let $g:\pi_1(W)\to G$ be
the epimorphism carrying the standard generators of $\pi_1(W)$ to
$\alpha_1, \beta_1, \dots, \alpha_d, \beta_d$, respectively.
Formula (\ref{Turaev1}) yields that
\begin{equation}\label{formula_Q_8_genus_d1}
|{\rm Hom}_g(\pi, G')| = 8^{2d-1}( 1 + \sum_{{i=1,2,3},\, G_{\chi_i}=G}
\, g^*(\zeta_{\chi_i}) ([W])\, +  \frac{1}{2^{2d-2}} \, g^*(\zeta_{\rho_0}) ([W]) ).
\end{equation}
Since  all values of the homomorphisms $\{\chi_i:Q_8\to
\CC^\times\}_{i=1,2,3}$  are $\pm 1$, the same is true for the
cocycles $\{\zeta_{\chi_i}\}_{i=1,2,3}$. Hence $g^*(\zeta_{\chi_i})
([W]) = \pm 1$ for all $i$. Then the sum $ 1 + \sum_i
g^*(\zeta_{\chi_i}) ([W])$ on the right-hand side of
(\ref{formula_Q_8_genus_d1}) is an integer. Since
$g^*(\zeta_{\rho_0}) ([W])$ is a root of unity and $d\geq 2$, the
number $|{\rm Hom}_g(\pi, G')|$ is non-zero.   The  same arguments
as in the proof of Theorem~\ref{Theorem_Q_8(1)} show that $$ |{\rm
Hom}_g(\pi, G')| = 8^{2d-1}( 1 + \left\{
                                                                             \begin{array}{ll}
                                                                               0  \\
                                                                               1   \\
                                                                               3
                                                                             \end{array}
                                                                           \right\} \pm \left(\frac{1}{2}\right)^{2d-2}
).$$\qed

\section{Proof of Theorems 3 and 4}

The group ${\mathcal A}={\Aut}(Q_8)$ is known to be isomorphic to
the symmetric group $S_4$. We shall  identify both ${\mathcal A}$
and $S_4$ with the group of rotations of a 3-dimensional cube as
follows. Let us label the vertices of the cube by $\{1,2,3,4\}$ such
that the vertices of each main diagonal have the same label, see
Figure~1. Let us label the faces of the cube   by $\pm i, \pm j, \pm
k$ so that   the opposite faces have opposite labels. Then any
rotation of the cube defines both a permutation in $S_4$ and an
automorphism of $Q_8$. This establishes the identification of these
groups mentioned above.

\begin{center}
\begin{picture}(160,130)%
\put(102,107){\bf 1} \put(2,107){\bf 2} \put(52,132){\bf 3}\put(152,132){\bf 4} 
\put(102,-5){\bf 3} \put(2,-5){\bf 4} \put(52,20){\bf 1}\put(152,20){\bf 2} 
\put(80,115){$i$} \put(130,70){{\it j}} \put(60,55){{\it k}} 
\put(80,15){$-i$} \put(25,70){{\it --j}} \put(90,75){{\it --k}} 
\put(5,5){\line(1,0){100}} \put(5,105){\line(1,0){100}} \put(55,130){\line(1,0){100}} 
\put(55,30){\line(1,0){45}} \put(110,30){\line(1,0){45}}
\put(5,5){\line(0,1){100}} \put(105,5){\line(0,1){100}} \put(155,30){\line(0,1){100}} 
\put(55,30){\line(0,1){70}} \put(55,110){\line(0,1){20}}
\qbezier(5,5)(15,10)(55,30) \qbezier(105,5)(115,10)(155,30) 
\qbezier(5,105)(15,110)(55,130) \qbezier(105,105)(115,110)(155,130) 
\put(65,-25){Figure 1}
\end{picture}
\end{center}

\bigskip

\bigskip

The stabilizers in ${\mathcal A}$ of the irreducible representations
$\{\chi_i\}_{i=0}^3$ and $ \rho_0$  of $Q_8$  can  be explicitly
computed. As mentioned above, ${\mathcal A}_{\chi_0}= {\mathcal
A}_{\rho_0}={\mathcal A}$.

\begin{lemma}\label{stabilizer} Let $\{{\mathcal
A}_{\chi_i}\subset \mathcal A\}_{i=1}^{3}$ be the stabilizers of the non-trivial
1-dimensional representations  $\{\chi_i\}_{i=1}^{3}$ of $Q_8$. Then
\begin{itemize}
\item[a.] Every  permutation of order $4$ in ${\mathcal
A}= S_4$ belongs to exactly one of the  groups $\{{\mathcal
A}_{\chi_i}\}_{i=1}^{3}$.
\item[b.] A permutation of order $3$ in ${\mathcal
A}=S_4$ belongs to none of the groups $\{{\mathcal
A}_{\chi_i}\}_{i=1}^{3}.$
\end{itemize}
\end{lemma}

{\it Proof.} The action of ${\mathcal A}$ on the set of non-trivial
1-dimensional representations of $ Q_8 $ is equivalent to the action on the kernels
 of these representations. These kernels  are  precisely the
order 4 subgroups $\langle i \rangle$, $\langle j \rangle$,  and $\langle k
\rangle$  of $Q_8$. The action of ${\mathcal A}$ on these subgroups  is
equivalent to the action of ${\mathcal A}$ on the pairs of opposite faces
$\{ \pm i\} , \{ \pm j\}, \{ \pm k\}$ of the cube.

Now,  a permutation of order 4 in $S_4$ corresponds to a rotation of the cube
about a line connecting  the centers of two opposite faces.
Such a rotation  stabilizes this pair of faces  and permutes the other  two
pairs. This implies  the first claim  of the lemma.

A permutation of order 3 in $S_4$ corresponds to a rotation of the cube about a main
diagonal. Such a rotation permutes cyclically
the pairs of opposite faces. This implies  the second  claim  of the lemma. \qed

{\it Proof of Theorem \ref{Theorem_Q_8(2)}.}  For each
$n=0,8,16,24,$ and 40 we produce two integers $k, m\geq 2$ and a
short exact sequence of groups
\begin{equation}\label{n=16}  Q_8 \hookrightarrow G'
\stackrel{q}{\twoheadrightarrow} \ZZ_k(\alpha) \times
\ZZ_m(\beta) \end{equation} such that $|L(\alpha,\beta,q)|=n$. Here
$\ZZ_k(\alpha)$ denotes the cyclic group $\langle \alpha \, | \,
\alpha^k = 1 \rangle$.

\bigskip

\noindent {\bf Case $n=40$.} Set $G = \ZZ_2(\alpha) \times
\ZZ_2(\beta)$ and $G'=Q_8\times G$. The map $q:G'\to G$ is the
projection.  Clearly,  $L(\alpha,\beta,q)= \{ \gamma_1, \gamma_2
\in Q_8\, |\, \gamma_1 \gamma_2 = \gamma_2 \gamma_1\}$. By Formula
(\ref{FroMed}), we have $|L(\alpha,\beta,q)|=40$.

\bigskip

For $n= 8,16,24$, we proceed as follows. We first produce a sequence
(\ref{n=16}) and  commuting lifts $\widetilde{\alpha},
\widetilde{\beta}\in G'$ of $\alpha$ and $\beta$. This allows us to
use Formula (\ref{formula_Q_8}) for  $|L(\alpha,\beta,q)|$. Set
$\Delta=\kappa_{\widetilde{\alpha}, \widetilde{\beta}}\in H_2(G')$.
The sign ${\omega}^*(\eta_{\rho_0}) (\Delta)=\pm 1$ in
(\ref{formula_Q_8}) will be   computed
 from the
equation
\begin{equation} \label{pairing}
{\mathcal M}_{\omega(\widetilde{\beta})} {\mathcal M}_{\omega(\widetilde{\alpha})}  =
{\omega}^*(\eta_{\rho_0}) (\Delta) {\mathcal M}_{\omega(\widetilde{\alpha})}
 {\mathcal M}_{\omega(\widetilde{\beta})}\, ,
\end{equation} where ${\mathcal M}_{\omega(\widetilde{\alpha})}, {\mathcal M}_{\omega(\widetilde{\beta})} \in GL_2(\CC)$ are
any matrices
  satisfying
\begin{equation} \label{matrix_Q_8}
\rho_0 (\widetilde{\alpha} \gamma \widetilde{\alpha}^{-1}) = {\mathcal M}_{\omega(\widetilde{\alpha})} \rho_0 (\gamma)
{\mathcal M}_{\omega(\widetilde{\alpha})}^{-1}, \mbox{ and } \rho_0 (\widetilde{\beta} \gamma \widetilde{\beta}^{-1}) = {\mathcal M}_{\omega(\widetilde{\beta})} \rho_0 (\gamma)
{\mathcal M}_{\omega(\widetilde{\beta})}^{-1},
\end{equation}
for all $\gamma \in Q_8$ (cf.\ (\ref{matrix})). Formula
(\ref{pairing}) follows   from the definition of $\eta_{\rho_0}$ and
the equalities  $
 \widetilde{\alpha}\widetilde{\beta}   =
 \widetilde{\beta}\widetilde{\alpha} $ and
$${\omega}^*(\eta_{\rho_0}) (\Delta) ={\omega}^*(\eta_{\rho_0}) (\kappa_{\widetilde{\alpha},
\widetilde{\beta}}) =
 \frac{ \eta_{\rho_0} (\omega(\widetilde{\alpha}), \omega(\widetilde{\beta}))}{ \eta_{\rho_0} (\omega(\widetilde{\beta}), \omega(
  \widetilde{\alpha}))}\, .$$

\bigskip

\noindent {\bf Case  $n=16$.} Set $G = \ZZ_3(\alpha) \times
\ZZ_3(\beta)$ and $G'=Q_8\rtimes (\ZZ_3(y) \times \ZZ_3(z))$, where
$y$ acts on $Q_8$ by $y(i)=j, y(j)=-k$, and $z$ acts trivially. The
map $q:G'\to G$ is given by $q(Q_8)=1$, $q(y)=\alpha$, and
$q(z)=\beta$. Clearly, $\widetilde{\alpha} = y$ and
$\widetilde{\beta} = z$ commute in $G'$.

Since $y$ acts as an automorphism of order 3 on $Q_8$, we have that
$\omega(G') \nsubseteq {\mathcal A}_{\chi_i}$ for $ i=1,2,3$ by
Lemma \ref{stabilizer}.  Since the action of $\widetilde{\beta}=z$
on $Q_8 $ is trivial, we have $ \omega(\widetilde{\beta})=1$ and
hence  ${\omega}^*(\eta_{\rho_0}) (\Delta) = 1$. Therefore
$|L(\alpha,\beta,q)| = 8(1 +1)=16$.

\bigskip

\noindent {\bf Case  $n=24$.} Set $G = \ZZ_2(\alpha) \times
\ZZ_2(\beta)$ and $G'=Q_8\rtimes (\ZZ_2(y) \times \ZZ_2(z))$, where
$\ZZ_2(y) \times \ZZ_2(z)$ acts on $Q_8$ by inner automorphisms:
\begin{equation}\label{inner}
y(\gamma)=j\gamma j^{-1}\,\,\,  \mbox{ and } \,\,\,
z(\gamma)=i\gamma i^{-1}\,\,\,  \mbox{ for all } \, \gamma \in
Q_8.
\end{equation}
The map $q:G'\to G$ is given by $q(Q_8)=1$, $q(y)=\alpha$,
and $q(z)=\beta$. Clearly, $\widetilde{\alpha} = y$ and
$\widetilde{\beta} = z$ commute in $G'$.

We have  $\omega(G') \subset  \Inn(Q_8) \subset {\mathcal
A}_{\chi_i}$ for $i=1,2,3$.  It follows from (\ref{inner}) that the
matrices
$${\mathcal M}_{\omega(y)} = \rho_0(j) = \left(
                                                            \begin{array}{cc}
                                                              0 & 1 \\
                                                              -1 & 0 \\
                                                            \end{array}
                                                          \right)
\quad \mbox{ and } \quad {\mathcal M}_{\omega(z)} = \rho_0(i)  = \left(
                                                            \begin{array}{cc}
                                                              i & 0 \\
                                                              0 & -i \\
                                                            \end{array}
                                                          \right)
$$ satisfy   (\ref{matrix_Q_8}).
By (\ref{pairing}),  ${\omega}^*(\eta_{\rho_0}) (\Delta) = -1$.
Thus, $|L(\alpha,\beta,q)|=8(1+3-1) = 24$.

\bigskip

\noindent {\bf Case  $n=8$.} Set $G = \ZZ_2(\alpha) \times
\ZZ_2(\beta)$ and $G' = Q_{16}\rtimes \ZZ_2(z)$, where $Q_{16}$ is
generated by $a,b$ subject to the relations $a^8 =1, a^4=b^2,
bab^{-1}=a^{-1}$, and   $z$ acts on $Q_{16}$   by $z(a) = a^3$ and
$z(b)=b$. We identify the  group $\langle a^2,b\rangle \subset G'$
with $Q_8$ via   $a^2=i$ and $b=j$. The map $q:G'\to G$ is given by
$q(a)=\alpha$, $q(b)= 1$,  and $q(z)=\beta$.  It is easy to check
that $\widetilde{\alpha}=ab$ and $\widetilde{\beta}= a^2bz$ commute
and  $q(\widetilde{\alpha})=\alpha$, $q(\widetilde{\beta})=\beta$.

Note that the conjugation by $a$ induces an 
automorphism of order 4 of $Q_8$, and $z$ acts by an inner
automorphism on $Q_8$:
$$ z(\gamma)=j\gamma j^{-1}\,\,\,  \mbox{ for all } \, \gamma
\in Q_8.$$
 It follows by Lemma
\ref{stabilizer} that $\omega(G')\subset {\mathcal A}_{\chi_i} $ for
exactly one $i\in \{1,2,3\}$. Since
$$\widetilde{\alpha} i \widetilde{\alpha}^{-1} = -i \quad \mbox{ and }
 \quad\widetilde{\alpha} j \widetilde{\alpha}^{-1} = k,$$
the matrix ${\mathcal M}_{\omega(\widetilde{\alpha})}$ has to
satisfy
$$\rho_0 (-i) =
{\mathcal M}_{\omega(\widetilde{\alpha})} \rho_0 (i) {\mathcal M}_{\omega(\widetilde{\alpha})}^{-1}
 \quad \mbox{ and } \quad  \rho_0 (k) =
{\mathcal M}_{\omega(\widetilde{\alpha})} \rho_0 (j) {\mathcal M}_{\omega(\widetilde{\alpha})}^{-1}.
$$
We take the following solution:  ${\mathcal
M}_{\omega(\widetilde{\alpha})} = \left(
                                                            \begin{array}{cc}
                                                              0 & 1 \\
                                                              i & 0 \\
                                                            \end{array}
                                                          \right).$
Similarly, since
$$\widetilde{\beta} i \widetilde{\beta}^{-1} = i  \quad \mbox{ and }  \quad \widetilde{\beta} j \widetilde{\beta}^{-1} = -j,$$
we can take ${\mathcal M}_{\omega(\widetilde{\beta})} = \left(
                                                            \begin{array}{cc}
                                                              1 & 0 \\
                                                              0 & -1 \\
                                                            \end{array}
                                                          \right).$
By (\ref{pairing}), ${\omega}^*(\eta_{\rho_0}) (\Delta)= -1$. Hence
$|L(\alpha,\beta,q)| = 8(1+1-1) = 8$.

\bigskip

\noindent {\bf Case  $n=0$.}    Set  $G'' = Q_{8}\rtimes \ZZ_6(y)$,
where   $y$ acts on $Q_{8}$   by $y(i) = j$ and $y(j)=-k$. Set $G' =
G'' \rtimes \ZZ_2(z) $, where   $z$ acts on $G''  $   by $z(i)=i$,
$z(j) = -j$ and $z(y) = k y$. Clearly, $Q_8$ is a normal subgroup of
$G'$  and $G'/Q_{8}=\ZZ_6 \times \ZZ_2$. We have a short exact
sequence
\begin{equation}\label{op}  Q_8 \hookrightarrow G'
\stackrel{q}{\twoheadrightarrow} G= \ZZ_6(\alpha) \times
\ZZ_2(\beta)\, , \end{equation} where  $q(Q_8)=1$,
$q(y)=\alpha$, and $q(z)=\beta$.


We use Formula (\ref{Turaev1-}) to prove that
$L(\alpha,\beta,q)=\emptyset$.  By Lemma \ref{stabilizer}, since the
conjugation by $y$ induces an  automorphism  of $Q_8$ of order 3, it
can  not stabilize  non-trivial 1-dimensional
representations of $Q_8$. It follows that only $\chi_0$ and $\rho_0$ contribute to
(\ref{Turaev1-}):
\begin{equation}  |L(\alpha,\beta,q)|=  8 \,
\left( g^*(\zeta_{\chi_0}) ([W]) + g^*(\zeta_{\rho_0})
([W])\right).
\end{equation}
As mentioned in Section~\ref{section2}, we have  $ g^*(\zeta_{\chi_0})
([W])=1$. Since $g^*(\zeta_{\rho_0}) ([W])$ is a root of unity and $|L(\alpha,\beta,q)|$ is an integer,  $
g^*(\zeta_{\rho_0}) ([W]) = \pm 1$. We compute  $g^*(\zeta_{\rho_0}) ([W]) $  from the definition of $\zeta_{\rho_0}$ in  Section~\ref{section2}. Put
$\widetilde{\alpha} = y$, $\widetilde{\beta} =  z$, and
$\widetilde{\alpha\beta} = \widetilde{\alpha}\widetilde{\beta}$.
Since
$$\widetilde{\alpha} i \widetilde{\alpha}^{-1} = j\ \mbox{\ \ and \ \ }\ \widetilde{\alpha} j \widetilde{\alpha}^{-1} = -k,$$
the matrix $M_{\alpha}\in GL_2(\CC)$ must satisfy (cf. (\ref{bh}))
$$M_{\alpha} \rho_0 (i) =  \rho_0 (j) M_{\alpha}
\ \mbox{\ \ and \ \ }\
M_{\alpha} \rho_0 (j) = \rho_0 (-k) M_{\alpha}.
$$
We can take $M_{\alpha} = \left(
                                                            \begin{array}{cc}
                                                              1 & -1 \\
                                                              i & i \\
                                                            \end{array}
                                                          \right).$
Note that
$\widetilde \beta \gamma \widetilde \beta^{-1}=z(\gamma)=i\gamma i^{-1}$ for all $\gamma \in Q_8$. We can take
$M_{\beta} = \left(
                                                            \begin{array}{cc}
                                                              1 & 0 \\
                                                              0 & -1 \\
                                                            \end{array}
                                                          \right)$.

The conjugation by $\widetilde{\alpha\beta} =
\widetilde{\alpha}\widetilde{\beta}$ induces the automorphism
$\displaystyle \left\{
                           \begin{array}{l}
                             i \mapsto j \\
                             j \mapsto k
                           \end{array}
                         \right.
$ of $Q_8$, and hence   $M_{\alpha\beta}$ must satisfy
$$M_{\alpha\beta} \rho_0 (i) =  \rho_0 (j) M_{\alpha\beta}
\ \mbox{\ \ and \ \ }\
M_{\alpha\beta} \rho_0 (j) = \rho_0 (k) M_{\alpha\beta}.
$$ We can take $M_{\alpha \beta} = \left(
                                                            \begin{array}{cc}
                                                              1 & 1 \\
                                                              i & -i \\
                                                            \end{array}
                                                          \right)$. 

We have $  (\widetilde{\alpha\beta})^{-1} \widetilde{\alpha}
\widetilde{\beta}  = 1$ and $ (\widetilde{\beta\alpha})^{-1}
\widetilde{\beta} \widetilde{\alpha} = z^{-1}y^{-1}zy = b = j
  \in Q_8$.  Thus, $\zeta_{\rho_0}(\alpha, \beta)$ and $ \zeta_{\rho_0}(\beta, \alpha)$ satisfy

$$
 \left(
\begin{array}{cc}
1 & 1 \\
i & -i \\
\end{array}\right)   = \zeta_{\rho_0}(\alpha, \beta) \left(
\begin{array}{cc}
1 & -1 \\
i & i \\
\end{array}\right)
 \left( \begin{array}{cc}
1 & 0 \\
0 & -1 \\
\end{array}\right)
$$ and
$$
 \left(
\begin{array}{cc}
1 & 1 \\
i & -i \\
\end{array}\right)   \left(
\begin{array}{cc}
0 & 1 \\
-1 & 0 \\
\end{array}\right) =  \zeta_{\rho_0}(\beta, \alpha)
\left( \begin{array}{cc}
1 & 0 \\
0 & -1 \\
\end{array}\right) \left(
\begin{array}{cc}
1 & -1 \\
i & i \\
\end{array}\right).$$
Therefore  $\zeta_{\rho_0}(\alpha, \beta)=1$, $\zeta_{\rho_0}(\beta,
\alpha)= -1$, and
$$g^*(\zeta_{\rho_0}) ([W]) =g^*(\kappa_{\alpha, \beta})= \frac{\zeta_{\rho_0}(\alpha, \beta)}{\zeta_{\rho_0}(\beta, \alpha)} = -1.$$
Thus $|L(\alpha,\beta,q)| = 8(1 -1)=0$.  \qed

\bigskip

\noindent {\bf Remark.} It is possible to obtain
$|L(\alpha,\beta,q)|= 24$ as a combination $  8(1+1+1)$ in Formula
(\ref{formula_Q_8}). Indeed, set $G = \ZZ_2(\alpha) \times
\ZZ_2(\beta)$ and $G' = Q_{16}\times \ZZ_2(z)$.   The map $q:G'\to
G$ is given by   $q(a)=\alpha$, $q(b)=1$, and $q(z)=\beta$. As in
Case $n=8$ above, $\omega(G')$ lies  in  exactly one of the
stabilizers $\{{\mathcal A}_{\chi_i}\}_{i=1}^{3}$. We choose
$\widetilde{\alpha} = a$ and $\widetilde{\beta} = z$. Since
conjugation by $\widetilde{\beta}=z$  is the identity on $Q_8=\Ker\,
q$, we have $\omega^*(\eta_{\rho_0})(\kappa_{\widetilde \alpha,
\widetilde \beta})= 1$. Therefore $|L(\alpha,\beta,q)|=8(1+1+1) $.

\bigskip

{\it Proof of Theorem \ref{Theorem_Q_8(2+)}.} Consider   the short
exact sequence (\ref{op}). As we know, $L(\alpha, \beta,
q)=\emptyset$. We claim that the induced homomorphism $q_* :
H_2(G'{}) \rightarrow H_2(G{})$ is an epimorphism so that
$\kappa_{\alpha, \beta}\in \im (q_*)$. Set
 $\Gamma=Q_8$ and consider the   exact sequence (see \cite{Br})
$$H_2(G')
\stackrel{q_*} \rightarrow H_2(G) \rightarrow \Gamma/[\Gamma,G']
\rightarrow H_1(G') \rightarrow H_1(G) \rightarrow 1\, .$$ The group
$ \Gamma/[\Gamma,G']$ is trivial because $[\Gamma,\Gamma] = \{ \pm
1\}$, and, as we saw, there is an element of $G'$ whose action on
$\Gamma$ cyclically permutes the non-trivial cosets of
$[\Gamma,\Gamma]$.  Therefore $q_*$ is an epimorphism.\qed

\section{Some other non-abelian groups}

In this section we state analogues of  Theorem \ref{Theorem_Q_8(1)}
for the symmetric group   $S_3$,   the dihedral group of order eight
$ D_4 $, the   alternating group of order twelve $A_4$, and the
 extra-special 2-groups.

\noindent {\bf 1.}
 Let $1 \rightarrow S_3 \rightarrow G'
\stackrel{q}\rightarrow G \rightarrow 1$ be a short exact sequence
of groups. For any   commuting elements $\alpha, \beta$ of  $G$, we
have  $|L(\alpha, \beta, q)|=18$. Indeed, the group $S_3$ has two
1-dimensional  representations $\chi_0$ and $\chi_1$, and
one 2-dimensional  irreducible representation $\rho_0$. Clearly, $G$ stabilizes them all.
Formula (\ref{Turaev1-}) gives
$$ |L(\alpha, \beta, q)| = |S_3|\, ( 1 +  g^*(\zeta_{\chi_1}) ([W])\, +
  g^*(\zeta_{\rho_0}) ([W]) ).
$$
Since $(S_3)_{ab} \cong \ZZ_2 $,   the values of $\chi_1$  are $\pm
1$. Hence $g^*(\zeta_{\chi_1}) ([W]) = \pm 1$, and, since
$g^*(\zeta_{\rho_0}) ([W])$ is a root of unity, $|L(\alpha, \beta,
q)|\neq 0$. We can therefore  apply Formula (\ref{Turaev 2}). Since
${\Aut}(S_3)=S_3$ and $H^2(S_3; \CC^\times) = 1$, the cohomology
class  $\eta_\rho$ arising from any irreducible representation
$\rho$ of $S_3$ is trivial. Hence $|L(\alpha, \beta,
q)|=6(1+1+1)=18$.

\medskip

\noindent {\bf 2.}  Let $1 \rightarrow D_4 \rightarrow G'
\stackrel{q}\rightarrow G \rightarrow 1$ be a short exact sequence
of groups and  $\alpha, \beta  $ be commuting elements of $G$. Then
$|L(\alpha, \beta, q)|=24$ or 40. Indeed, $\Aut (D_4)=\ZZ_2 \times
\ZZ_2$, and all automorphisms of $D_4$ are inner. Hence the action
of $G$ on $\Irr (D_4)$ is trivial. The group $D_4$ has four
1-dimensional  representations $\{\chi_i\}_{i=0}^3$ and
one 2-dimensional  irreducible representation
 $\rho_0$.   By   (\ref{Turaev1-}),
$$|L(\alpha, \beta, q)| = |D_4|\, ( 1 + \sum_{i=1,2,3}
\, g^*(\zeta_{\chi_i}) ([W])\, + \, g^*(\zeta_{\rho_0}) ([W]) ),$$
where $g^*(\zeta_{\chi_i}) ([W])=\pm 1$ and $g^*(\zeta_{\rho_0})$ is
a root of unity. Consequently, $L(\alpha, \beta, q) \neq \emptyset$
and we can apply Formula (\ref{Turaev 2}), where
${\omega}^*(\eta_{\chi_i}) (\Delta)=1$ for all $i=1,2,3,$ and
${\omega}^*(\eta_{\rho_0}) (\Delta)=\pm 1$. Therefore $|L(\alpha,
\beta, q)|=8(4 \pm 1)=24$ or 40.


\medskip

\noindent {\bf 3.}   Let $1 \rightarrow A_4 \rightarrow G'
\stackrel{q}\rightarrow G \rightarrow 1$ be a short exact sequence
of groups and  $\alpha, \beta  $ be commuting elements of $G$.  Then
$\vert L(\alpha, \beta, q)\vert=0,24,$ or 48.  Indeed, the group
$A_4$ has three 1-dimensional  representations and a
3-dimensional irreducible representation  $\rho_0$. Clearly, if
${\Aut}(A_4)=S_4$ stabilizes a non-trivial 1-dimensional
representation of $A_4$, then it stabilizes all of them.
 If $\vert L(\alpha, \beta, q)\vert \neq 0$, then  by (\ref{Turaev 2}), we have  $|L(\alpha, \beta, q)|=12(1+
\left\{
                                                                             \begin{array}{l}
                                                                               0 \\
                                                                               2
                                                                             \end{array}
                                                                           \right\} \pm 1)
 $ so that $|L(\alpha, \beta, q)|=24 $ or 48.

\bigskip Theorem \ref{Theorem_Q_8(1)} and Example 2 above in this section can be
generalized as  follows. Recall that  a finite  group $\Gamma$ is
extra-special if its order is a power of a prime integer $p$,  its  center $Z(\Gamma)$  is isomorphic to $\ZZ_p$,
 and the quotient $\Gamma/Z(\Gamma)$ is a direct product of several copies of $\ZZ_p$. An
extra-special $p$-group $\Gamma$ satisfies $[\Gamma, \Gamma]=Z(\Gamma)$ and   $\vert \Gamma\vert=p^{2r+1}$ for some $r$, see
\cite[Theorem 5.2]{G}. For example, both $Q_8$ and $D_4$ are
 extra-special 2-groups.

\medskip

\noindent {\bf 4.} Let $\Gamma$ be an extra-special 2-group of order
$2^{2r+1}$ with $r\geq 1$. Then $\Gamma$ has $| \, \Gamma/[\Gamma,\Gamma] \, | =
2^{2r}$ one-dimensional  representations. Since $Z(\Gamma)=\ZZ_2$, the group $\Gamma$ has two
conjugacy classes or order 1. All other  conjugacy classes in $\Gamma$ have  exactly two elements.
Indeed, if  $Z(\Gamma) = \{1,\sigma\} $ with $\sigma\in \Gamma$, then the conjugacy class of $\gamma \in G-
Z(\Gamma)$  is equal to   $
\{\gamma, \sigma\gamma \}$. It follows that $\Gamma$ has $2^{2r}+1$
conjugacy classes. Hence $\Gamma$  has a unique irreducible
representation $\rho_0$ of dimension greater than 1. The equality
$\sum_{\rho\in \Irr\Gamma} (\dim \rho)^2 =|\Gamma|$ implies that
$\dim \rho_0=2^r$.

Let $1 \rightarrow \Gamma \rightarrow G' \stackrel{q}\rightarrow G
\rightarrow 1$ be a short exact sequence of groups and  $\alpha,
\beta  $ be commuting elements of $G$.    Formula (\ref{Turaev 2}) implies that either $\vert L(\alpha, \beta,
q)\vert = 0$ or
$$|L(\alpha, \beta, q)| = 2^{2r+1}(1+ N \pm 1)
$$
where $N$ is the number of non-trivial
1-dimensional  representations of $\Gamma$ stable under the action of
$G$ by conjugations.


\begin{thebibliography}{x}
%

\bibitem{Br} K. Brown, Cohomology of groups, \newblock{\it Graduate Texts in
Mathematics, 87}, Springer-Verlag, New York-Berlin, 1982.

\bibitem{FQ} G. Freed, F. Quinn, Chern-Simons theory with finite gauge
group,
\newblock{\it Comm. Math.  Phys.,}  {\bf 156} (1993), 435--472.

\bibitem{Fr} G. Frobenius, \"{U}ber Gruppencharaktere,
\newblock{\it Sitzber. K\"{o}niglich Preuss. Akad. Wiss. Berlin,}
(1896), 985--1021.


\bibitem{G} D. Gorenstein,
\newblock Finite Groups,
\newblock Chelsea, New York, 1980.

\bibitem{Jo} G. Jones, Enumeration of homomorphisms and surface-coverings,
\newblock{\it Quart. J. Math. Oxford Ser. (2),} {\bf 46}  (1995),  no. 184, 485--507.


\bibitem{Me} A. Mednykh, Determination of the number of nonequivalent
coverings over a compact Riemann surface.
\newblock{\it Dokl. Akad. Nauk SSSR,} {\bf 239} (1978),   269--271 (Russian).
\newblock{\it Soviet Math.\ Dokl.,} {\bf 19} (1978), 318--320 (English
transl.).

\bibitem{T} V. Turaev, On certain enumeration problems in two-dimensional topology, arXiv:0804.1489

\end{thebibliography}
\end{document}